\documentclass[10pt,notitlepage,a4paper,reqno,ps]{amsart}

\usepackage[utf8]{inputenc}
\usepackage[T1]{fontenc}
\usepackage[english]{babel}
\usepackage{amsmath}
\usepackage{amsfonts}
\usepackage{amssymb}
\usepackage{amsthm}
\usepackage{bbm}
\usepackage{mathrsfs}
\usepackage{graphicx}
\usepackage{braket}
\usepackage{color}
\usepackage{amscd,epsfig,esint,graphics}
\usepackage{wrapfig}
\usepackage{verbatim}
\usepackage{mathtools}
\usepackage{color}

\theoremstyle{plain}
\newtheorem{theorem}{Theorem}[section]
\newtheorem{lemma}[theorem]{Lemma}

\newtheorem{remark}[theorem]{Remark}
\newtheorem{definition}[theorem]{Definition}

\theoremstyle{definition}
\theoremstyle{remark}
\numberwithin{equation}{section}

\newcommand{\R}{\mathbb{R}}

\newcommand{\N}{\mathbb{N}}

\newcommand{\om}{\Omega}

\DeclareMathOperator*{\esssup}{ess\,sup}

\definecolor{vg}{rgb}{0.0, 0.26, 0.15}

\title{Relaxation of one-dimensional nonlocal supremal functionals in the Sobolev setting} %and dimension reduction}
\author[A.\,Torricelli]
{Andrea Torricelli}
\address[A.\,Torricelli]{Dipartimento di Scienze Fisiche, Informatiche e Matematiche, Universit\`a degli Studi di Modena e Reggio Emilia, via Campi 213/b, 41125, Modena, Italy.}
\email{andrea.torricelli@unipr.it}

\author[E.\,Zappale]
{Elvira Zappale}
\address[E.\,Zappale]{Dipartimento di Scienze di Base ed Applicate per l'Ingegneria, Sapienza-Universit\`a di Roma, via A. Scarpa, 16, 00161 Roma, Italy}
\email{elvira.zappale@uniroma1.it}

\begin{document}
	
	\baselineskip 3.4ex
	\vspace{0.5cm}
	\maketitle
	
	\begin{abstract}
		
		We provide necessary and sufficient conditions on the density $W:\R^d\times\R^d\to\R$ in order to ensure the sequential weak* lower semicontinuity of the functional $J: W^{1,\infty}(I;\mathbb R^d)\to \mathbb R$, defined as \begin{align*}
			J(u):=\esssup_{I\times I}W(u'(x), u'(y)),
		\end{align*}
		when $I$ is an open and bounded interval of $\mathbb R$.
		We also show that, when $d=1$, the lower semicontinuous envelope of $I$ in general can be
		obtained by replacing $W$ by its separately level convex envelope.  
		

	\end{abstract}

%%%%%%%%%%%%%%%%%%%%%%%%%%%%%%%%%%%%%%%%%%%%%%%

\section{Introduction}

 A growing interest is developping towards nonlocal functionals in recent years, due to their many applications, both from the application view point (peridynamics, image processing, artificial intelligence, etc), and the theoretical one, e.g. \cite{BMCP, MD, DFKS, EDLL14, KS, CKS} among a pletora of scientific contributions.
 In particular a big attention is devoted to nonlocal integrals, we refer to \cite{Pm} and the bibliography contained therein for an overview and to \cite{BDM2} for the interactions with local energies, among a wide literature.
 
 On the other hand, $L^\infty$ variational problems arise when studying  optimal design problems, in the control theory, optimal transport, etc, and in these past decades a wide theory has been built, we refer to \cite{BL, ABPrin} for the pioneering theoretical papers and to \cite{EP, PZ, CK}  and the bibliography contained therein just to mention some more theoretical recent contribution among a much wider literature.    
 Our contribution inserts in the framework of detecting necessary and sufficient conditions for the lower semicontinuity and relaxation of nonlocal $L^\infty$ functionals depending on the derivatives of Sobolev fields defined on the real line.

In this paper we  focus on $L^\infty$ variational models, for which one might be interested in detecting equilibrium configurations, which consist of determining the lower semicontinuity of suitable functional, or, when this is not the case in their relaxation. When $\Omega$ is a bounded open subset of $\mathbb R^n$, $n > 1$, it is  well known that in general there is no explicit representation for the relaxation of nonlocal supremal functionals of the type
\begin{equation}\label{tildeJ}\tilde J(u):=\esssup_{\om\times \om}W(\nabla u(x), \nabla u(y))
	\end{equation}
with respect to the $W^{1,\infty}-weak^*$ topology. As anticipated, the question has been completely answered in \cite{KRZ} in the case where $\nabla u$ is replaced by $u$, while, in the case of gradients, a sufficient condition on $W$ ensuring lower semicontinuity has been provided in \cite{GZ}.

The current analysis, focusing on fields defined on an interval strongly rely on the results contained in \cite{KRZ}. 
% very much related with the analogous one treated in \cite{KRZ}. DISCUSS MORE ANG MAKE COMMENTS AND REFERENCES.
The study of the lower semicontinuity of $\tilde J$ with respect to the  $W^{1,\infty}-weak^*$ topology is  equivalent to characterize, for every $c \in \mathbb R$, the $W^{1,\infty}-weak^*$ closure of $B_{L_c(W)}$, with  $L_c(W)$ being $c$-(sub)level set of $W$, i.e.
\begin{align*}
	L_c(W):=\left\{(\xi, \eta) \in \R^{d\times n}\times\R^{d\times n}: W(\xi,\eta)\le c\right\},
\end{align*}
and, for every $K \subset \R^{d \times n} \times \R^{d\times n}$
\begin{align}\label{AK}
	B_K:=\left\{u \in W^{1,\infty}(\om;\mathbb R^d): (\nabla u (x), \nabla u(y)) \in K \text{ for a.e. }x\in \om\times\om \right\}.
\end{align}
In this note, stemming from the results in \cite{KZ, KRZ}, we provide an answer when $n=1$ in terms of the notion of Cartesian convexity introduced in \cite{KRZ} for the set $K$ and an equivalent characterization in term of separate convexity when also $d=1$. 
%\color{red}
%Indeed a function is weakly lower semicontinuous if and only if all its level sets are weakly closed. It is easy to see that $B_{L_c(W)}$ is the $c$-level set of $J$, so our result can be seen as a necessary and sufficient condition on $W$ that ensure that $B_{L_c(W)}$ is $W^{1,\infty}-weak^*$ closed. \\

We will also provide the supremal counterpart of the results for nonlocal integral functionals depending on derivatives of Sobolev functions contained in \cite{BP}. In that paper the authors find that the necessary and sufficient condition for the $L^p$ weak lower semicontinuity of functionals such as
\begin{equation}
	\label{nonlocalint}
	\iint_{I\times I}f(u'(x),u'(y))dxdy, \quad u\in W^{1,p}(I),\, I \text{ interval of }\R,
\end{equation}
is the separate convexity of $f$.

We will characterize the $W^{1,\infty}$-lower semicontinuity of the functional $\tilde J$ in \eqref{tildeJ} in terms of properties on the supremand $W$ when the open set $\Omega$ is an open interval on the real line.  For the sake of exposition we will mainly focus on the case when $d=1$, i.e. the field $u$ is scalar valued. On the other hand the results concerning fields in $W^{1,\infty}(I;\mathbb R^d)$ can be deduced exploiting analogous arguments, see Theorem \ref{cnsvec}. Actually, Theorem \ref{cns} could be deduced as a corollary of Theorem \ref{cnsvec}, but we opted for presenting the proof in the scalar valued setting for the sake of exposition.

We emphasize also that, while in \cite{BP} it has been proven that the relaxation of a nonlocal integral functional defined in open subsets of $\mathbb R$, does not admit, in general, a representation of the same type (see also \cite{BM-C, KZdintegral, BDM3}), we will show that, under suitable assumptions on $W$, in the case $n=1$, the relaxed functional,  can be represented in supremal form with a suitable density. Furthermore if also $d=1$, without any restriction on $W$, the relaxed energy is of the same type with a density coinciding with the separately level convex envelope of the original one. In the vector valued case the relaxation will be deduced in view of \cite[Theorem 1.2]{KRZ}, see Theorem \ref{relaxvec}.

To ease the parallel with the nonlocal integral setting, it is worth to mention that in \cite{M}, Mu$\tilde{n}$oz proved that separate convexity of the integrand is a necessary and sufficient condition for the weak lower semicontinuity of the functional \eqref{nonlocalint} when it depends on $\nabla u$, with $u\in W^{1,p}(\Omega)$ and $\Omega\subset \R^n$. The same question in the supremal setting is currently still open and it is not clear if the Cartesian level convexity of the supremand $W$ (see Section 3 for the definition introduced in \cite{KRZ}) is also necessary for the lower semicontinuity of $\tilde J$. 
We also stress that, in contrast with the integral framework, even in the local scalar setting, i.e. $n=1$ or $d=1$, the convexity notions arising in the supremal framework are not equivalent, see \cite{RZ2}, indeed this latter aspect prevents us from using the same strategies adopted in \cite{M}.

%On the other hand, in \cite{KZdintegral} the authors proved that the relaxation of such nonlocal integral functional is not structure preserving, i.e. its relaxation is not a nonlocal integral functional.

From now on, we will denote with $I$ the interval $(a,b)$, for some $a,b\in\R$, $a<b$. Referring to Section \ref{pre} for definitions and properties of the nonlocal densities below, 
the first result that we prove is the following 
\begin{theorem}
	\label{cns}
	Let  $J_{1d}:W^{1,\infty}(I)\to\R$ be defined as
	\begin{align}
		\label{functional1d}
		J_{1d}(u):=\esssup_{I\times I}V(u'(x), u'(y)),
	\end{align} with $V:\mathbb R \times \mathbb R \to \mathbb R$ coercive, lower semicontinuous, diagonal and symmetric, then $J_{1d}$ is lower semicontinuous with respect to the $W^{1,\infty}-weak^*$ topology if and only if $V$ is separately level convex.
\end{theorem}
Moreover the following result holds: 
\begin{theorem}
	\label{relax}
	Let $V$ and $J_{1d}$ be as in Theorem \ref{cns}. %defined as in \eqref{functional1d} with $V$ coercive, lower semicontinuous, diagonal and symmetric.
	Let $J_{1d}^{rlx}$ be the sequentially $W^{1,\infty}-$weak$^*$ lower semicontinuous envelope of $J_{1d}$, i.e.
	\begin{equation}\label{Jrlx1d}
	J_{1d}^{rlx}:=\inf\{\liminf_{j\to +\infty} J_{1d}(u_j): u_j \overset{*}{\rightharpoonup} u \hbox{ in }W^{1,\infty}(I)\}.
	\end{equation}
	Then
	$$J^{rlx}_{1d}(u)=\esssup_{I\times I} V^{slc}(u'(x), u'(y)),$$ for every $u \in W^{1,\infty}(I)$, where $V^{slc}$ is the separately level convex envelope of $V$.
\end{theorem}
For what concerns the vector valued case, we have
\begin{theorem}
	\label{cnsvec}
	Let  $J:W^{1,\infty}(I;\R^d)\to\R$ be defined as
	\begin{align}
		\label{functionald}
		J(u):=\esssup_{I\times I}W(u'(x), u'(y)),
	\end{align} with $W:\mathbb R^d \times \mathbb R^d \to \mathbb R$ coercive, lower semicontinuous, diagonal and symmetric, then $J$ is lower semicontinuous with respect to the $W^{1,\infty}-weak^*$ topology if and only if $W$ is Cartesian level convex.
\end{theorem}

Furthermore 
\begin{theorem}
	\label{relaxvec}
	Let $W$ and $J$ be as in Theorem \ref{cnsvec}. %defined as in \eqref{functional1d} with $V$ coercive, lower semicontinuous, diagonal and symmetric.
	Let $J^{rlx}$ be the sequentially $W^{1,\infty}-$weak$^*$ lower semicontinuous envelope of $J$, i.e.
%	\begin{equation}\label{Jrlxd}
$		J^{rlx}:=\inf\{\liminf_{j\to +\infty} J(u_j): u_j \overset{*}{\rightharpoonup} u \hbox{ in }W^{1,\infty}(I;\R^d)\}. $
%	\end{equation}
	If every sublevel set of $W$ has a basic Cartesian convexification, then
	$$J^{rlx}(u)=\esssup_{I\times I} {W}^{\times lc}(u'(x), u'(y)),$$ for every $u \in W^{1,\infty}(I;\mathbb R^d)$, where ${W}^{\times lc}$ stands for the Cartesian level convex envelope of $ W$.
\end{theorem}

	As already mentioned above, our results can be understood as results on non-local differential inclusions. Indeed, let $E$ be a compact subset of $\mathbb R\times \mathbb R$, then Theorem \ref{relax} states that for every sequence $(u_j)_j\subset W^{1,\infty}(I)$ such that $u_j\overset{*}{\rightharpoonup}u$ in $W^{1,\infty}(I)$ and
	\begin{align*}
		(u'_j(t),u'_j(s))\in E
	\end{align*}
for a.e. $(s,t) \in I\times I$, then 
\begin{align*}
	(u'(t),u'(s))\in {\widehat E}^{sc},
\end{align*}
where
\begin{align}\label{Ehat}
	%	&K^{sym}:= \left\{ (\xi, \eta)\in K: (\eta, \xi) \in K \right\},\\ %&K^{diag}:=\left\{(\xi,\eta)\in K: (\xi,\xi),(\eta,\eta)\in K\right\},\\ &
	\widehat{E}:=\left\{(\xi,\eta)\in E: (\eta,\xi),(\eta,\eta),(\xi, \xi)\in E\right\}%=K^{sym}\cap K^{diag},	
\end{align}
and the supersript {\rm sc} stands for its separately level convex hull.

Consequently Theorem \ref{cns} states that 
	$(u'(t),u'(s))\in E$ for a.e. $(s,t)\in I\times I$ if and only if $E$ is separately convex.

The differential inclusions counterpart of Theorems \ref{cnsvec} and \ref{relaxvec} will be discussed in Section \ref{proofs}. 

The paper is organized as follows: in section \ref{pre} we start fixing notation and  providing some preliminary results. The proofs of the main results are given in the last section.

\section{Preliminaries and Notation}\label{pre}
In the sequel $n$ ad $d$ will denote elements of $\mathbb N$. In this section we will give some defintion and recall some known results that will be useful for the proofs of our theorems.
 
\color{black}

Moreover we recall the notions of \textit{level convexity, separate convexity} and \textit{separate level convexity.}
\begin{definition}
	\label{sep_conv}
	A set $A \subset \R^d\times\R^d$ is said  {\it separately convex} if for every $(\xi_1,\xi_2), (\eta_1,\eta_2) \in \R^d\times \R^d$ with $\xi_1=\eta_1$ or $\xi_2=\eta_2$ it holds
	\begin{align*}
		t(\xi_1,\xi_2)+(1-t)(\eta_1,\eta_2) \in A,
	\end{align*}
	for every $t \in (0,1).$ We denote with $A^{sc}$ the separate convex hull of $A$, namely the smallest separately convex set containing $A$.
\end{definition}

\begin{remark}
	\label{compactness}
	If $A\subset \mathbb R^d\times \mathbb R^d$ is open then so is $A^{sc}$. This is in general not true for compactness, see \cite[Remark 7.18 (ii)]{Dac}, but, as observed in \cite[Proposition 2.3]{K} the property is preserved if $A\subset \R\times\R$. 
\end{remark}
For any function $f:\R^d\times\R^d \to \R \cup \{\infty\}$ and any real number $c$, by $L_c(f)$ we denote the level set $c$ of $f$, i.e.
$$
L_c(f):=\{(\xi, \eta) \in \R^d \times \R^d: f(\xi, \eta)\leq c\}.
$$
\begin{definition}
	\label{sep_level_conv}
	A function $f:\R^d\times\R^d \to \R \cup \{\infty\}$ is said to be \begin{itemize}
		\item [1.] {\it level convex} if for every $c \in \R$ the set $L_c(f)$ is convex; 
		\item[2.]  {\it separately convex} if for every $\xi \in \R^d$ both $f(\cdot, \xi)$ and $f(\xi, \cdot)$ are convex functions. 
		\item[3.]  {\it separately level convex} if for every $c \in \R$  $L_c(f)$ is separately convex.
	\end{itemize}
	We denote the separately  (level) convex envelope of $f$, namely the greatest (level) convex function less than or equal to $f$ by ($f^{slc}$) $f^{sc}$.
\end{definition}

\begin{definition}\label{WdiagDEF}
	Given $W:\R^d\times\R^d\to\R\cup \{\infty\}$ we say that $W$ is symmetric if $W(\xi,\eta)=W(\eta, \xi)$ for every $(\xi,\eta)\in\R^d\times\R^d$. Also, we say that $W$ is diagonal if 
	\begin{align*}
		(\xi,\eta)\in L_c(W) \Rightarrow (\xi,\xi),(\eta,\eta)\in L_c(W).
	\end{align*}
	In other words, following \cite[eq. (28)]{KRZ},
	\begin{align*}
		W(\xi,\eta)=\max\{W(\xi,\xi), W(\xi,\eta), W(\eta, \eta)\}.
	\end{align*}
\end{definition}

\begin{definition}
	We say that a function $W:\R^d\times\R^d\to\R\cup\{\infty\}$ is coercive if there exists $C'>0$ such that
	\begin{align*}%\label{coerciW}
		C'|(\xi,\eta)| \leq W(\xi,\eta),
	\end{align*}
	for every $\xi,\eta\in\R^d.$
\end{definition}
\begin{remark}
	\label{identity}
	From Definitions \ref{sep_conv} and \ref{sep_level_conv},  for any $W:\R^d\times\R^d\to\R$,
	\begin{align*}
		L_c(W)^{sc}\subset L_c(W^{slc}),
	\end{align*}
	The inverse inclusion in general does not hold, while if $d=1$, in  \cite[Lemma 7.4]{KZ} it has been proven that if $W$ is symmetric and diagonal then 
	\begin{align*}
		L_c(W)^{sc}=L_c(W^{slc}).
	\end{align*}
\end{remark}

\subsection{Properties of the supremal nonlocal functionals}

In this subsection we will describe the properties of the energy densities appearing in Theorems \ref{cns}- \ref{relaxvec}.
First of all we observe that there is no loss of generality in supposing that the supremands $V$ and $W$ in \eqref{functional1d} and \eqref{functionald} are diagonal and symmetric.

In fact, let $W: \mathbb R^{d\times n} \times \mathbb R^{d\times n} \to \mathbb R$ and  define the function $\widehat W:\mathbb R^{d\times n} \times \mathbb R^{d\times n} \to \mathbb R$ by
\begin{align}
	\label{simdiag}
	\widehat{W}(\xi,\eta):=\inf\left\{ c\in\R: (\xi,\eta)\in\widehat{L_c(W)} \right\},
\end{align}
as in \cite[(7.1)]{KZ}, with $\widehat{L_c(W)}$ given by \eqref{Ehat}. $\widehat{W}$ is the symmetric and diagonal, and $\widehat{W}\ge W, $ so the coercivity of $W$ is inherited by $\widehat{W}$. 
Following \cite{KZ}, define for $\Omega \subset \mathbb R^n$, and every $K \subset \R^{d\times n} \times \R^{d\times n}$
\begin{align}\label{AKdef}
	A_K:=\left\{v \in L^{\infty}(\om;\mathbb R^{d\times n}): (v (x), v(y)) \in K \text{ for a.e. }x\in \om\times\om \right\}.
\end{align} 
In \cite[Proposition 5.1]{KZ} it has been proven that $A_K = A_{\widehat K}$, from which follows 
\cite[eq (7.3)]{KZ}, i.e. for every $u \in W^{1,\infty}(\Omega;\mathbb R^d)$,
%by \cite[eq (7.2), Proposition 5.1 and eq. (7.4)]{KZ} ,
\begin{align*}
	\esssup_{\om\times\om}\widehat{W}(\nabla u(x),\nabla u(y))=%\inf\left\{ c\in\R: \nabla u\in A_{L_c(\widehat{W})} \right\}\\
	%&=\inf\left\{ c\in\R: \nabla u\in A_{\widehat{L_c(W)}} \right\}=
	%\inf\left\{ c\in\R: \nabla u\in A_{L_c(W)} \right\}\\
	\esssup_{\om\times\om}W(\nabla u(x),\nabla u(y)).
\end{align*}
In particular the last equality holds when $W=V:\mathbb R\times \mathbb R\to\R$ and $\Omega=I$.

\color{black}

%This shows that the functional $J$ is invariant by symmetrization and diagonalization of its supremand.

%	On another note, for a generic $W$ symmetrical and diagonal we have that
%	\begin{align*}
%		L_c(W)=\widehat{L_c(W)},
%	\end{align*}
%	where we recall that $\widehat{L_c(W)}$ is the diagonalization and symmetrization of $L_c(W)$. 
%\end{remark}

\color{black}
\section{Proofs}\label{proofs}
In this section, first we prove that the separate level convexity of the supremand $V$ appearing in \eqref{functional1d} %being separately level convex 
is a sufficient and necessary condition for the sequential $W^{1,\infty}-$weak$^*$ lower semicontinuity of fhe functional $J_{1d}$  defined in \eqref{functional1d}.
%$J_{1d}(u):=\esssup_{(x,y) \in I\times I}V(u'(x), u'(y))$.
To this end, we recall  for every $E \subset  \mathbb R\times \mathbb R$, the set \begin{align}\label{BEdef}
	B_E:=\left\{u \in W^{1,\infty}(I): (u' (x), u'(y)) \in E \text{ for a.e. }(x,y)\in I\times I \right\}.
\end{align} 
Moreover, we recall that $J_{1d}$ is sequentially $W^{1,\infty}-$ weakly$^*$ lower semicontinuous if and only if $B_{L_c(V)}$ is sequentially $W^{1,\infty}-$weakly$^*$ closed for every $c \in \mathbb R$. With the following results we will prove that $B_{L_c(V)}$ is sequentially $W^{1,\infty}-$weakly$^*$ closed if and only if $B_{L_c(V)}=B_{L_c(V)^{sc}}$ and by Remark \ref{identity} we have that this is true if and only if $V$ is separately level convex.

The next lemma is key to prove of Theorem \ref{cns}. Let $I=(a,b)\subset\R$, we will denote with $S^\infty_E(I)$ the set of all the simple functions $s$ (see for instance \cite[eq (5.6)]{KZ}), and with $D^\infty(I)$ the set of all the integral functions
	
	\begin{align}\label{integralfunction}
		u_s(x,c):=c+\int_a^x s(t)dt,\quad x\in I,
	\end{align}
	with $s\in S^\infty(I)$ and $c\in\R$. 
\begin{remark}
	\label{Amerio}
	We recall that,  given $v\in L^\infty(I)$, the function $u_v(\cdot,c) \in L^\infty(I)$ defined, for every $c\in\R$ as \eqref{integralfunction}, i.e.
	$$
	u_v(x,c):=c+ \int_a^x v(t)dt,
	$$  
	is an element of $W^{1,\infty}(I)$. Moreover, the distributional derivative exists also in the classical sense and
	\begin{align*}
		\frac{du_v(x,c)}{dx}=v(x)\quad \text{for a.e. }x\in I.
	\end{align*}
Since for our subsequent analysis the role of $c$ can be neglected, we will omit it from the notation and denote the integral function just by $u_v$.

Thus, given a function $u_v$ with $v\in L^\infty(I)$ and a sequence of functions $(v_j)_{j\in\N}\subset L^\infty(I)$ such that $v_j\overset{*}{\rightharpoonup} v$ in $L^\infty(I)$, then $u_{v_j}\overset{*}{\rightharpoonup}u_v$ in $W^{1,\infty}(I)$ (see \cite{F}).  In particular if $v$ is a simple function $s\in S^\infty(I)$, it is immediately verified that $D^\infty(I)\subset 
W^{1,\infty}(I).$

\end{remark}

\begin{lemma}
	\label{dense}
	Given $E\subset\R\times\R$ symmetric and diagonal, let $B_E$ be as in \eqref{BEdef}, then, for every $u\in B_E$ there exists a sequence $(u_j)_{j\in\N}\subset B_E\cap D^\infty(I)$ such that $u_j\overset{*}{\rightharpoonup}u$ in $W^{1,\infty}(I)$.
\end{lemma}
\begin{proof}
	Since $B_E\subset W^{1,\infty}(I)$ we can 
define $B'_E:=\{u': u \in  B_E\}$. Clearly $B'_E\subset L^\infty(I)$ and it is a subset of 
$A_E:=\{v \in L^\infty(I): (v(x), v(y))\in E \}$ in \eqref{AKdef}.
Hence, by \cite[Lemma 5.4 and Corollary 5.5]{KZ} any element  $v \in A_E$ can be approximated in $L^\infty(I)$ strong convergence by a sequence $\{s_j\}_j \subset A_E\cap S^\infty(I)$. In particular this happens for $v=u' \in B'_E$.
Thus, by Remark \ref{Amerio}, given $u \in B_E$, there exists $\{u_{s_j}\}_j\subset D^\infty_I\cap B_E$
weakly* converging to $u$ in $W^{1,\infty}(I)$.
\end{proof}

\begin{remark}
In the previous lemma we supposed that $E \subset \R \times \R$ is symmetric and diagonal, but from \cite[Proposition 5.1]{KZ} we can deduce that the previous lemma holds for every $E$ closed and contained in $\mathbb R^{d\times n} \times \mathbb R^{d\times n}$, with $n, d \in \mathbb N$.
\end{remark}
\color{black}

\begin{proof}[Proof of Theorem \ref{cns}]
	Due to the coercivity assumption on $V$ and the metrizability of the weak* topology on bounded sets we will omit the word 'sequentially'.
	The sufficiency of the separate level convexity of $V$ follows from one of the implications of \cite[Theorem 1.3(i)]{KZ}, taking into account the weak* convergence of $(u'_k)_k$ for every sequence $(u_k)_k$ weakly* converging to $u$ in $W^{1,\infty}(I)$.
	
	For what concerns the reverse implication, we claim \color{black} that the $W^{1,\infty}-$weak* lower semicontinuity of the functional $J_{1d}$ in \eqref{functional1d}  ensures the $L^\infty-$weak* lower semicontinuity of the functional in $L^\infty(I)$ defined as $\esssup_{I\times I}
V(v(x), v(y))$ 	along sequences $(v_k)_k \subset L^\infty(I)$. Indeed, if $v_k \overset{\ast}{\rightharpoonup} v$, in $L^\infty (I)$, then, as observed in Remark \ref{Amerio}, it is possible to define \begin{align*}%\label{uk}
	u_k(x):=\int_a^x v_k(t)dt \in L^\infty(I)
\end{align*} and $u_k$ is in $L^\infty(I)$. Indeed the fact that $v_k$ is bounded in $L^\infty(I)$, entails also that
$$
\esssup_{x \in I}|u_k(x)| \leq \esssup_{x \in I}\int_a^x |v_k|(x)dx \leq \|v_k\|_{L^{\infty}(I)}(b-a).$$ Moreover $\frac{d u_k}{dx} $ exists for a.e. $x \in I$ and coincides with $v_k$. 
Hence $u_k \in W^{1,\infty}(I)$. Moreover $u_k$ is unifomrly bounded in the $W^{1,\infty}$ norm, thus, up to a not relabelled subsequence, $u_k \overset{\ast}{\rightharpoonup} u$ in $W^{1,\infty}(I)$ with $\frac{d u}{dx}= v$ a.e. in $I$. This proves our claim.
%The above observations allow us to assert the weak* lower semicontinuity of  $\esssup_{I\times I}
%V(v(x), v(y))$. 
Consequently the proof is concluded by the reverse implication of \cite[(i)Theorem 1.3]{KZ}.   
 \end{proof}
\color{black}

\begin{remark}\label{remvec}
It is worth to observe that the level convexity of $W:\mathbb R^{d\times n}\times \R^{d\times n}\to \mathbb R$ is sufficient to ensure the weak* lower semicontinuity of the functional $\esssup_{\om\times \om}W(\nabla u(x),\nabla u(y))$ in view of the same arguments exploited in the first part of the proof of Theorem \ref{cns}.  In view of \cite[Theorem 1.1]{KRZ} the reverse implication does not hold, even if one assumes $n=1$: we refer to Theorem \ref{cnsvec} below.
\end{remark}

Now we prove that the relaxed functional $J_{1d}^{rlx}$ in \eqref{Jrlx1d} is given by
\begin{align*}
	%\label{thm17thesis}
	J_{1d}^{rlx}(u)=\esssup_{ I\times I}V^{slc}(u'(x),u'(y)).
\end{align*}

\begin{proof}[Proof of Theorem \ref{relax}]
	
The proof follows by \cite[Proposition 7.5]{KZ} and Theorem \ref{cns}. We present some details for the reader's convenience. 
Let us define 
\begin{align*}
	J^{slc}(u):=\esssup_{I\times I}V^{slc}(u'(x),u'(y)).
\end{align*}
By definition of separately level convex envelope (see Definition \ref{sep_level_conv}) we have that 
\begin{align*}
	V\ge V^{slc},
\end{align*}
so, by \cite[Theorem 1.3]{KZ} for every $(u_k)_{k\in\N}\subset W^{1,\infty}(I)$ such that $u_k\overset{*}{\rightharpoonup}u$ in $W^{1,\infty}(I)$ we have that $u'_k \overset{\ast}{\rightharpoonup}u'$ in $L^\infty(I)$, hence
\begin{align*}%\label{lb}
	J^{slc}(u)\le\liminf_{k \to \infty}J^{slc}(u_k)\le\liminf_{k \to \infty}J_{1d}(u_k).
\end{align*}
To prove the other inequality, fix $u\in W^{1,\infty}(I)$. 

%we can define
%\begin{align*}
%	c:=J^{slc}(u)<+\infty,
%\end{align*}
%and we can choose a sequence $(c_k)_{k\in\N}\subset\R$ such that $c_k \searrow c$ as $k\to\infty$ and such that
%\begin{align*}
%	u\in B_{L_c{_k(V^{slc})}}=B_{L_c{_k(V)^{sc}}},
%\end{align*}
%where in the last equality we have exploited \cite[Lemma 7.4]{KZ}.

By \cite[Theorem 1.3]{KZ} %Proposition \ref{cns} 
there exists a sequence $(v_k)_{k\in\N}\subset L^{\infty}(I)$ such that $v_k\overset{*}{\rightharpoonup}u'$ in $L^{\infty}(I)$ as $k\to+\infty$ and
\begin{align}\label{almostrecovery}
\lim_{k \to +\infty}\esssup_{I\times I} V(v_k(x),v_k(y))= \esssup_{I\times I} V^{slc}(u'(x), u'(y)).
\end{align}
Again the possibility of defining  $u_k(x):=\int_a^x v_k(x)dx$, see Remark \ref{Amerio}, allows us to say that $(u_k)_{k\in \mathbb N} \subset W^{1,\infty}(I)$, and $u_k \overset{\ast}{\rightharpoonup} w$ in $W^{1,\infty}(I)$ with $\frac{d w}{dx}= u'$, a.e. in $I$, hence $w(x)= u (x)+ C$. This fact, together with \eqref{almostrecovery} guarantee that
$$
\lim_{k\to +\infty} \esssup_{I\times I}V(u_k'(x), u'_k(y))=\esssup_{I\times I}V^{slc}(u'(x), u'(y)),
$$ 
which ensures the existence of a recovery sequence, i.e. $(u_k)_{k \in \mathbb N} \subset W^{1,\infty}(I)$, weakly* converging to $u$ and that concludes the proof.
\end{proof}

In \cite{KRZ} the following notions have been introduced.
\begin{definition}\label{defvec}
	A set $E \subset \R^d \times \R^d$ is {\it Cartesian convex} if $A \times A \subset E$ implies $A^{co} \times A^{co} \subset  E$, where $A^{co}$ denotes the convex envelope of the set $A$.
	
	If $E$ is not Cartesian convex, the smallest Cartesian convex set containing $E$ is, by definition, its Cartesian convex hull and it is denoted by $E^{\times sc}$. 
	
	$E$ admits a basic Cartesian convexification if 
	$E^{\times xc}= \bigcup_{A \times A \in \mathcal P_E} A^{co}\times A^{co}$,
	where $\mathcal P_E$ represents the set of all the maximal Cartesian squares i.e. sets $B\times B$, contained in $E$ and with the property that if there exists a set $C$ such that $C\times C \subset E$ and $B\subset C$, then $B=C.$ 
	
	A function $f:\R^d\times\R^d \to \R \cup \{\infty\}$ is said to be 
	{\it Cartesian level convex} if for every $c \in \R$ the set $L_c(f)$ is Cartesian convex.

	When $f$ is not Cartesian level convex, the corresponding envelope, namely 
	$$f^{\times lc}(\xi, \eta) := \sup\{g (\xi, \eta) |\, f : \R^d \times \R^d \to \R \hbox{ is Cartesian level convex and }g \leq f\},
	$$ 
can be described as
	$$f^{\times lc}(\xi, \eta) = \inf\{c \in \mathbb R : (\xi, \eta) \in L_c(f)^{\times c}
	\}$$ for $(\xi, \eta) \in \R^d \times \R^d$.
\end{definition}

\begin{proof}[Proof of Theorem \ref{cnsvec}]
	The proof follows along the lines of Theorem \ref{cns}, considering vector valued functions in $W^{1,\infty}(I;\mathbb R^d)$, replacing \cite[Theorem 1.3(ii)]{KZ} by \cite[Theorem 1.1]{KRZ}, i.e. exploiting the necessary and sufficient condition given by the Cartesian level convexity of $W$ for nonlocal supremal functionals depending on $L^\infty(I;\mathbb R^d)$ fields, instead of the separate level convexity.   
\end{proof}

We also observe, as already mentioned in the introduction, that Theorem \ref{cns} is a corollary of Theorem \ref{cnsvec}, since Cartesian level convexity reduces to separate level convexity when $d=1$, as proven in \cite{KRZ}.

In terms of differential inclusions, given a compact set $K\subset\R^d\times\R^d$ we denote the $W^{1,\infty}-$weak$^*$ closure of $B_K$ in \eqref{AK} as $B^\infty_K$, i.e.
\begin{align}\label{Bkinfty}
	B^\infty_K:=\left\{ u\in W^{1,\infty}(I;\mathbb R^d)| \exists (u_j)_{j\in\N}\subset B_K : u_j \overset{*}{\rightharpoonup}u \in W^{1,\infty}(I;\mathbb R^d) \right\}.
\end{align}

Theorem \ref{cnsvec} states that $B^\infty_K= B_K$ if and only if $K$ is Cartesian convex.

The same arguments of Theorem \ref{relax} can be applied to describe the relaxation of the functional \ref{functionald} in the case when $W$ admits a {\it basic Cartesian convexification}

\begin{proof}[Proof of Theorem \ref{relaxvec}]
	The proof relies on the same arguments of Theorem \ref{relax}, exploiting \cite[Theorem 1.2]{KRZ} instead of \cite[Theorem 1.3]{KZ}. 
	\end{proof}

\begin{remark}
	\label{lastrem}
	It is worth to observe that when the diagonality assumption on the supremands $V$ and $W$ above is dropped the relaxation results hold in terms of the envelopes of $\widehat V$ and $\widehat W$, respectively, defined accordingly to \eqref{simdiag}.
\end{remark}

Clearly Theorem \ref{relaxvec} in terms of differential inclusions says that $B^\infty_K$  in \eqref{Bkinfty} coincides with $B_{K^{\times c}}$ if $K$ is symmetric, diagonal and admitting a basic Cartesian convexification.
 
We conclude by observing that all the above results could be rephrased in terms of the lower semicontinuity or relaxation of nonlocal indicator functionals of the type $$\iint_{I\times I}\chi_K (u'(x), u'(y)) dx dy,$$
where for $K \subset \R^d \times \R^d$, $\chi_K$ stands for its characteristic function defined as 
$	\chi_K(\xi,\eta):=
	\begin{cases}
		&0 \hbox{ if }(\xi, \eta)\in K,\\
		&\infty \hbox{ otherwise}.
	\end{cases} $
This leads to a characterization of the lower semicontinuity for the above unbounded nonlocal integrals in terms of the Cartesian convexity of $K$ and a characterization of its relaxation in terms of its Cartesian convex envelope, if and only if $K$ satisfies the assumptions of \cite[Theorem 1.2]{KRZ}. 
 
\bigskip

{\bf Acknowledgements}
The authors are members of INdAM-GNAMPA and they gratefully acknowledge the support of Progetto 'GNAMPA Prospettive nelle scienze dei materiali: modelli variazionali, analisi asintotica e omogeneizzazione'.
EZ thanks Dipartimento di Scienze Fisiche, Informatiche e Matematiche of University of Modena and Reggio Emilia for its support and kind hospitality.

\color{black}

\end{document}